\documentclass[12pt]{article}
\usepackage[english]{babel}
\usepackage{authblk}
\usepackage[hmarginratio=3:2]{geometry}
\usepackage{graphicx}
\usepackage{amsmath}
\usepackage{amsthm}
\usepackage{amssymb}
\usepackage[mathscr]{eucal}
\usepackage{textcomp}
\usepackage[latin1]{inputenc}
\usepackage{latexsym}
\usepackage{hyperref}
\usepackage[dvipsnames]{xcolor}
\usepackage{framed}
\usepackage{tikz-cd}
\usepackage{ifsym}
\newtheorem{theorem}{Theorem}[section]
\newtheorem{lemma}[theorem]{Lemma}
\newtheorem{corollary}[theorem]{Corollary}
\newtheorem{proposition}[theorem]{Proposition}

\title{On the Chow ring of the classifying stack of algebraic tori}
\author[1]{Francesco Sala}
\affil[1]{\textit{Scuola Normale Superiore, Piazza dei cavalieri 7, 56126, Pisa, Italy}}
\date{}
\begin{document}
\maketitle

\begin{abstract}
We investigate the structure of the Chow ring of the classifying stacks $BT$ of algebraic tori, as it has been defined by B. Totaro.\\
Some previous work of N. Karpenko, A. Merkurjev, S. Blinstein and F. Scavia has shed some light on the structure of such rings. In particular Karpenko showed the absence of torsion classes in the case of permutation tori, while Merkurjev and Blinstein described in a very effective way the second Chow group $A^2(BT)$ in the general case. Building on this work, Scavia exhibited  an example where $A^2(BT)_\text{tors}\neq 0$.\\
Here, by making use of a very elementary approach, we extend the result of Karpenko to special tori and we completely determine the Chow ring $A^*(BT)$ when $T$ is an algebraic torus admitting a resolution with special tori $0\rightarrow T\rightarrow Q\rightarrow P$. In particular we show that there can be torsion in the Chow ring of such tori.
\end{abstract}
\bigskip
\bigskip
\section*{Introduction}

In this paper we give some new examples of Chow rings of classifying spaces of
algebraic tori with non-trivial additive torsion; our example is more elementary than the first one due to F. Scavia \cite{Sca}, and it lies in $A^3(BT)$; the main tool is a description of the Chow group of the total space for a torsor under a permutation torus, which is of independent interest.

Let $G$ be an affine algebraic group over a field $k$. B. Totaro \cite{To} gave a definition of the Chow ring $A^*_G$ of the classifying space $BG$; when $k = \mathbb C$ there is a natural ring homomorphism from $A^*_G$ to the cohomology of $BG$, which is an isomorphism after tensoring with $\mathbb Q$, but not integrally (\cite{EG}, \cite{MV}).

The ring $A^*_G$ has been computed for many classes of reductive split groups, such as $GL_n$, $SL_n$, $SO_n$, $PGL_p$ when $p$ is a prime (see \cite{MV},\cite{Vi}). On the other hand, not much is known about Chow rings of classifying spaces of non-inner forms of these groups.

When $T = \mathbb G^n_m$ is a totally split algebraic torus, the Chow ring $A^*_
T$ is a polynomial ring $\mathbb Z[t_1, \dots , t_n]$, where $t_i$ is the first Chern class of the $i$-th projection $\mathbb G^n_m \rightarrow   \mathbb G_m$,
considered as a representation of $T$. In this paper we consider non-totally split algebraic tori; let $T$ be such a torus.

Choose a finite Galois extension $k'$ of $k$, with Galois group $\Gamma$, such that $T_{k'}\simeq\mathbb G^n_{m,k'}$. The group $\Gamma$ acts on $A^*_{T_{k'}}\simeq\mathbb  Z[t_1, \dots, t_n]$; the natural pullback map $A^*_T\rightarrow A^*_{T_{k'}}$ factors through the invariant subring $(A^*_{T_{k'}})^\Gamma$. We call the resulting homomorphism of graded rings $A^*_T\rightarrow  (A^*_{T_{k'}})^\Gamma$ \emph{the base-change homomorphism}. It is not hard to check that it is independent of the choice of $k'$, and that it is an isomorphism when tensored with $\mathbb Q$. It is a natural question whether it is always an isomorphism.

N. Karpenko \cite{Kar} showed that this is always the case when $T$ is a permutation torus. 

A. Merkurjev and S. Blinstein (see \cite{BiMe} and \cite{Mer}) gave a cohomological interpretation of the kernel $(A^2_T)_\text{tors}$ of the base-change homomorphism in degree $2$. By using this, Scavia gave an example of a torus $T$ for which $(A^2_T)_\text{tors} = 0$, so that the base-change homomorphism is not always injective in degree $2$.\\

Recall that an algebraic group $G$ is said to be \emph{special} if every $G-$torsor is Zariski-locally trivial; it is easy to see that permutation tori are special. In this paper we show that the base-change is an isomorphism for all special tori; then we investigate the simplest class of non-special algebraic tori, those tori $T$ that can be embedded into a special torus $P$, in such a way that the quotient $P/T$ is also special. For these tori Merkurjev \cite{Mer} showed that $(A^2_T)_\text{tors}$ is always $0$.

Here we give an example of tori of this special form for which the base-change
homomorphism is not injective nor surjective.\\

In the first section we study the simplest representation of a permutation torus $T$, showing the existence of an equivariant decomposition in $N=\dim(T)$ strata. These strata will be seen to be disjoint unions of permutation tori, defined over some suitable extension $k\subseteq L\subseteq k'$.

In the second section we show a general result (Theorem $2.1$) which allows to calculate the Chow ring $A^*(E)$ for a $T-$torsor $E\rightarrow M$ from the ring $A^*(M)$, valid for every permutation torus $T$. We show that the pull-back map $A^*(M)\rightarrow A^*(E)$ is surjective, and we give a description of its kernel; this is our main new contribution.

In the third section we briefly recall the results which are known from the literature and extend Karpenko's result to special tori. We finally concentrate on the $T$'s with a resolution by special tori $0\rightarrow T\rightarrow P\rightarrow Q\rightarrow 0$. We give a formula for the Chow ring $A^*_T$ as a function of $A^*_P,A^*_Q$, based on Theorem $2.1$.

Using this formula, in the fourth and last section we find explicit examples where the map $\phi:A^*_T\rightarrow (A^*_{T_{k'}})^\Gamma$ in not injective or not surjective.\\

This work arose from a suggestion that my advisor Angelo Vistoli gave me for my master thesis project. I would like to thank him gratefully for his guidance and support.

\section{Stratifying a representation of permutation tori}

Here we begin to set up the tools which we will need to calculate the Chow ring of the classifying stack $BT$ for a class of algebraic tori $T$. 

In particular we start by considering quasi-split tori over a fixed field $k$, i.e. algebraic tori which arise as Weil restriction from a finite type étale $k-$algebra. If $E$ is such an algebra, we thus denote $T:=R_{E/k}\mathbb{G}_{m,E}$. If $k'$ is a Galois closure of $E$, then $T_{k'}$ is a split torus, and the group $\Gamma:=\text{Gal}(k':k)$ acts by permutation on the character lattice of $T_{k'}$ (or, in other words, this lattice has a $\Gamma-$invariant basis on which the action of $\Gamma$ permutes the elements). Therefore, these groups are also classified under the name of \emph{permutation tori}.

Indeed let $E=\prod_{i=1}^m k_i$, where each $k_i$ is a finite separable extension of $k$ of degree $d_i$ (in particular $k'$ will contain the normal closure over $k$ of each $k_i$). If $T_i:=R_{k_i/k}\mathbb{G}_{m,k_i}$ then $T=\prod_{i=1}^m T_i$.

We now prove that each factor $R_{k_i/k}\mathbb{G}_{m,k_i}$ is isomorphic - as a scheme - to $\mathbb{A}^{d_i}_k-\{N^{k_i}_k=0\}$. Here $\mathbb{A}^{d_i}_k$ is $\text{Spec } k[x_1,\dots, x_{d_i}]$ and $N^{k_i}_k$ is the polynomial in the $x_j$'s given by $\prod_{l=1}^{d_i} \sigma_l(\sum \alpha_j x_j)$ - where $\alpha_1,\dots,\alpha_{d_i}$ is a basis of $k_i$ over $k$ and the $\sigma_l:k_i\hookrightarrow k'$ are the immersions of $k_i$ in $k'$ leaving $k$ fixed. This is known from basic field theory to have coefficients in $k$.

In order to prove the above claim let us fix an extension $k':k$, let $F$ be the Galois closure of $k'$ over $k$ and let $J$ be the set of immersions $k'\hookrightarrow F$ fixing $k$. Then $\Gamma:=\text{Gal}(F:k)$ acts on $J$ by permuting its elements, and we have a well-known 

\begin{proposition}
Let $X'$ be a variety over $k'$ and let $\overline{X}:=\underset{j\in J}{\prod}X'\times_{k',j}F$. For every $\sigma\in \Gamma$ let $\phi_\sigma$ be the automorphism of $\overline{X}$ induced by the action of $\sigma$ on $J$.

Then the $\phi_\sigma$ define Galois descent data, and the variety over $k$ defined by these data is $R_{k'/k}X'$.
\end{proposition}



With this result we can easily construct $R_{k'/k}X'$ for a given $k'-$variety $X'$. In particular let us concentrate on $X'=\mathbb G_{m,k'}$. 

Let $x_1,\dots, x_n$ be a set of indeterminates indexed by $J$, and let $e_1,\dots,e_n$ be a basis of $k'$ over $k$. Then $\underset{\sigma:k'\hookrightarrow F}{\prod}(\sum_{j\in J}\sigma(e_j)x_j)$  is the homogeneous polynomial $N^{k'}_k$ in the $x_j'$s which defines the norm over $k$ of the element $\sum x_je_j$. Let $X=\mathbb A^n_k-\{N^{k'}_k=0\}$.

In order to prove our claim, we can just consider $X\times_kF$ (with $F$ a Galois closure as above): this is equal to $\mathbb A^n_{F}-\{\prod_{\sigma:k'\hookrightarrow F}\sum e_jx_j=0\}\simeq \mathbb A^n_{F}-\{y_1\dots y_n=0\}$ via the isomorphism sending $\sum \sigma(e_i)x_i$ to $y_\sigma$. But it is exactly $\underset{j\in J}{\prod}\mathbb G_{m,k'}\times_{k',j}F$, with the action of $\Gamma=\text{Gal}(F:k)$ which permutes the factors as it does on $J$ (corresponding to the action of $\Gamma$ on the factors $\sum\sigma(e_i)x_i$). This completes the proof.\\

We note that each $R_{k_i/k}\mathbb G_{m,k_i}$ splits over $\overline k$ as a product of $[k_i:k]$ factors isomorphic to $\mathbb G_{m,\overline k}$, and $\Gamma$ permutes them according to its action on the set of $k-$immersions $k_i\hookrightarrow \overline k$.\\

Now, $R_{k_i/k}\mathbb{G}_{m,k_i}$ has a natural representation $V_i$  over $k$ of degree $d_i$ and the whole product $T$ has a natural representation $V:=\oplus V_i$ of degree $N:=\sum_id_i$. 

Then $V_{k'}:=V\times_{\text{Spec }k}\text{Spec } k'$ is a representation of the split torus $T_{k'}$ whose action is induced by the tautological inclusion $T_{k'}\hookrightarrow GL(V_{k'})$ of diagonal matrices.\\

From now on we will call $I_i$ the set of integers in the interval $(\sum_{j<i}d_j,\sum_{j\leq i}d_j]$, which can label the immersions $k_i\hookrightarrow\overline{k}$.\\

We have $T_i\times_{\text{Spec }k} \text{Spec }k'\simeq \mathbb{A}^{d_i}_{k'}-\{\prod_{l=1}^{d_i} \sigma_l(\sum \alpha_j x_j)=0\}$, which gives an isomorphism $T_i\times_{\text{Spec }k} \text{Spec }k'\simeq \mathbb{G}_{m,k'}^{ d_i}=\text{Spec }k'[y_1^{\pm 1},\dots y_{d_i}^{\pm 1}]$ sending $y_l$ to $\sigma_l(\sum \alpha_j x_j)$ for each $l$. In particular   $V_i\times_{\text{Spec }k} \text{Spec }k'\simeq \text{Spec }k'[y_1,\dots, y_{d_i}]$ inherits an action of $\Gamma$ which is the obvious one on $\text{Spec }k'$ and sends $y_l$ to $y_{\tau(l)}$, where $\tau$ is the permutation of $\{1,\dots, d_i\}$ induced by the action of $\Gamma$ on the $\sigma_l$'s.

To summarize, the group $\Gamma$ acts on the space $V\times_{\text{Spec }k} \text{Spec }k'$ as follows: tautologically on $\text{Spec }k'$, while for each $i$ it permutes $\{x_j : j\in I_i\}$ according to the identification of the latter with the set of $k-$immersions $k_i\hookrightarrow \overline k$.\\

Clearly $V=V_{k'}/\Gamma$, where the action of $\Gamma$ on $V$ is the one described at the end of the preceding paragraph.

For each $j\in\{1,\dots N\}$ let $H_j$ be the hyperplane in $V_{k'}=\text{Spec }k'[y_1,\dots , y_N]$ defined by $\{y_j=0\}$ and for each $J\subseteq \{1,\dots N\}$ let $H_J$ be the linear subspace of $V_{k'}$ obtained by $\underset{j\in J}{\bigcap} H_j$.

Then if $\overline{Z_p}:=\underset{|J|=p}{\coprod}H_J$, these closed subschemes are $T_{k'}$ and $\Gamma$-stable; moreover these two actions are compatible with the action of $\Gamma$ on $T_{k'}$, in the sense that for $g\in\Gamma, t\in T_{k'}$ it holds $g(t\cdot x)=g(t)\cdot g(x)$ for every $x$. This implies that 
the closed subschemes $Z_p:=\overline{Z_p}/\Gamma\hookrightarrow V$ are $T-$stable. Thus we have found a $T-$equivariant stratification of $V$, and the latter is the disjoint union of the $W_p:=Z_{p}-Z_{p+1}$, which we call \emph{the strata}. 

In particular $\overline{Z_0}=T_{k'}$ and $Z_0=T$.\\

We now describe quite explicitly the strata $W_{p}$, in terms of the action of $\Gamma$ on $\{1,\dots, N\}$.\\

Let us take $J\subseteq\{1,\dots, N\}$ such that $|J|=p$ and consider the split $N-p$-torus $T_J:=H_J-\underset{J\subsetneq I}{\bigcup} H_I$. Clearly $\Gamma$ acts on the set $\{I\subseteq \{1,\dots, N\}:|I|=p\}$ and the orbit of $J$ for this action $\text{Orb}_\Gamma(J)$ determines the orbit of $T_J$ for the action of $\Gamma$ on $\overline{Z_p}$: $\text{Orb}_\Gamma(T_J)=\underset{I\in\text{Orb}_\Gamma(J)}{\coprod}T_I.$

We shall see, by Galois descent, how to characterize $\text{Orb}_\Gamma(T_J)/\Gamma$; we will show that this is again a permutation torus defined over a finite extension $k_J:k$.\\

Let us consider for a moment a field $\mathfrak{k}$ such that $\mathfrak{k}\subseteq k_i$ for all $i$, and let $\mathfrak{\Gamma}\subseteq \Gamma$ be the group $\text{Gal}(k':\mathfrak{k})$ . Then $\Gamma$ acts on $R_{E/\mathfrak{k}}\mathbb{G}_{m,E}\times_{\text{Spec }k}\text{Spec }k'$, which is a disjoint union of $n:=[\mathfrak{k}:k]$ tori of dimension $N/n$, as follows: the orbits of $\mathfrak{\Gamma}$ partition $\{1,\dots, N\}$ into $n$ pieces, and we can identify each of these $n$ tori with one of them (where characters corresponds to elements $i\in\{1,\dots, N\}$); then $\Gamma$, acting on $\{1,\dots, N\}$ induces the sought action on these tori.

 This variety embeds $\Gamma-$equivariantly into $V_{k'}$, as we can see immediately. Indeed if $\mathfrak{I}\subset \{1,\dots,N\}$ is the subset corresponding to the immersions $k_i\hookrightarrow k'$ (for any $i$) which leave $\mathfrak{k}$ fixed, and $\mathfrak{I}^c$ is its complement, then $\text{Orb}_\Gamma(T_{\mathfrak{I}^c})$ is a disjoint union of $n$ split tori of dimension $N/n$ (since $\Gamma$ partitions $\{1,\dots,N\}$ into $N/n$ subsets isomorphic to $\mathfrak{I}$); in particular it is isomorphic to $R_{E/\mathfrak{k}}\mathbb{G}_{m,E}\times_{\text{Spec }k}\text{Spec }k'$. According to the observations we made in the preceding paragraphs, $\Gamma$ acts on this orbit exactly as it does on $R_{E/\mathfrak{k}}\mathbb{G}_{m,E}\times_{\text{Spec }k}\text{Spec }k'$, so that the embedding of the latter into $V_{k'}$ is $\Gamma-$equivariant.

In particular by Galois descent we have that $(\text{Orb}_\Gamma(T_{\mathfrak{I}^c}))/\Gamma\simeq R_{E/\mathfrak{k}}\mathbb{G}_{m,E}$. Moreover we note that the stabilizer of $\mathfrak{I}$ in $\Gamma$ is exactly $\mathfrak{\Gamma}$, while the $k_i$'s correspond to the subsets of $\mathfrak{\Gamma}$ which stabilize the orbits of its action on $\mathfrak{I}$ (or on $\mathfrak{I}^c$, which does not change them). Thus we can determine $\mathfrak{k}$ and all the $k_i$'s just by looking at the action of $\Gamma$ on $\text{Orb}_\Gamma(T_{\mathfrak{I}^c})$.\\

The same reasoning allows us to determine $(\text{Orb}_\Gamma T_J)/\Gamma$ for a general $J$. 

First of all, let $S_J\subseteq \Gamma$ be the stabilizer of $J$, which corresponds to a field extension $K:k$. Moreover let $J_1,\dots,J_s$ be the orbits of the elements of $J$ under the action of $S_J$, and $S_1,\dots,S_s$ the respective stabilizers. The $S_i$'s induce extensions $K_i:K$, and let $\mathbb{E}$ be the étale $K-$algebra $\prod_{i=1}^sK_i$. Then we can compare the action of $\Gamma$ on $R_{\mathbb{E}/K}\mathbb{G}_{m,\mathbb{E}}\times_{\text{Spec }k}\text{Spec }k'$ with the action of $\Gamma$ on $\text{Orb}_\Gamma(J)$, so that we are exactly in the situation of the last paragraph (we just restrict to the vector subspace of $V_{k'}$ generated by $\{y_i:i\in\underset{g\in\Gamma}{\cup} g(J)\}$).

We conclude that $\text{Orb}_\Gamma(T_J)/\Gamma\simeq R_{\mathbb{E}/K}\mathbb{G}_{m,\mathbb{E}}$. Let this torus be called $T_{\text{Orb}(J)}$.\\

So each of the $W_p$'s is a disjoint union of algebraic tori defined over extensions $k\subseteq K\subseteq k'$: $Z_p=\underset{|J|=p}{\coprod} T_{\text{Orb}(J)}$.\\

We conclude this section with the crucial

\begin{lemma} Let $J\subseteq\{1,\dots, N\}$, and $\mathfrak{k}$ be the fixed field inside $k'$ of the stabilizer in $\Gamma$ of $J$ and $J^c$.
Then $T_{\text{Orb}(J)}\times_{\text{Spec }\mathfrak{k}}T_{\text{Orb}(J^c)}=T\times_{\text{Spec }k}\text{Spec }\mathfrak{k}$.
\end{lemma}

\begin{proof} Let $\mathfrak{\Gamma}$ be the group $\text{Gal}(k':\mathfrak{k})$. As before, we can take the product (over $\text{Spec }\mathfrak{k}$) of both sides with $\text{Spec }k'$. Then $T\times_{\text{Spec }k}\text{Spec }\mathfrak{k}\simeq T_{k'}/\mathfrak{\Gamma}$, where the action of $\mathfrak{\Gamma}$ on $T_{k'}$ is just the one induced by the action of $\Gamma$. So by Galois descent we just need to verify that the action of $\mathfrak{\Gamma}$ on $(\text{Orb}_\Gamma(T_J)/\Gamma\times_{\text{Spec }\mathfrak{k}}\text{Orb}_\Gamma(T_{J^c})/\Gamma)\times_{\text{Spec }\mathfrak{\Gamma}}\text{Spec }k'$ (which is a split $N-$torus) can be identified with the action of $\mathfrak{\Gamma}$ on $T_{k'}$. But this is clear since the action of $\mathfrak{\Gamma}$ on $\text{Orb}_\Gamma(T_{J})/\Gamma\times_{\text{Spec }\mathfrak{\Gamma}}\text{Spec }k'$ is induced by the restriction of the action of $\Gamma$ on $T_{k'}$ to the subspace generated by $\{y_i:i\in J^c\}$; so the product of the $J,J^c$-tori gives a split $N-$torus with the entire action that $\mathfrak{\Gamma}$ had on $T_{k'}$.
\end{proof}

\section{The Chow ring of $T-$torsors}

We now exploit the results of the previous section to give a general formula for the Chow ring of a $T-$torsor, where $T$ is a permutation torus.

We will keep the same notation used before. In particular, $V$ will denote the standard representation of $T$, and the $Z_p$'s will be the closed sets of the stratification we found. Moreover we have a set of permutation tori $W_p=Z_p-Z_{p+1}$ defined over field extensions of $k$.\\
From now on, we will also call $U_p:=V-Z_p$ the open sets of the stratification.

Let $E\xrightarrow{p} M$ be a $T-$torsor of smooth schemes. Then we have a ring map $p^*:A^*(M)\rightarrow A^*(E)$ given by pull-back.

Clearly if $x\in A^*(M)$ is a characteristic class relative to $E$ (or, more fancifully, an element in the image of the pull-back $A^*(BT)\rightarrow A^*(M)$ for the map which classifies $E$) we have $p^*(x)=0$; indeed $p^*(x)$ is nothing but a characteristic class relative to the pull-back torsor $p^*E\rightarrow E$, which is the trivial torsor $E\times T\rightarrow E$.

Analogously, let us fix a field extension $L\supseteq k$, and let $E_L\xrightarrow{p_L} M_L$ be the corresponding $T_L-$torsor. Then, as before, we have that if $x\in A^*(M_L)$ is a $T_L-$characteristic class relative to $E_L$ then $p_L^*(x)=0$. Consider the commuting square
\[
\begin{tikzcd}
A^*(M_L)\arrow{r}\arrow{d}{p^*_L}&A^*(M)\arrow{d}{p^*}\\
A^*(E_L)\arrow{r}&A^*(E)
\end{tikzcd}
\]
where the horizontal arrows are the push-forwards for the obvious projections. Then, if $x'\in A^*(M_L)$ lies in the ideal generated by the $T_L-$characteristic classes induced by $E_L$ and if $x$ is its image in $A^*(M)$, we see immediately that $p^*(x)=0$ in $A^*(E)$.

We claim that the classes we found so far are the only generators of $\ker(p^*)$:

\begin{theorem}
Let $E\xrightarrow{p} M$ be a $T-$torsor. Let $\mathcal{I}\subseteq A^*(M)$ be the ideal generated by the push-forwards $\pi^L_*(x\cdot m)$, where for some field $L\supseteq k$ we denote $\pi^L:M_L\rightarrow M$ the projection and $x, m\in A^*(M_L)$, $x$ being a $T_L$-characteristic class induced by $E_L$.

Then $A^*(E)=A^*(M)/\mathcal{I}$. 
\end{theorem}

\begin{proof}

We will prove by descending induction on $p$ that $A^*(M)\rightarrow A^*(E\times_TU_p)$ is surjective, and that its kernel lies in $\mathcal{I}$. This will be sufficient, since $U_0=T$ and $E\times_TU_0=E$.

The base case is clear, since $U_N=V-\{0\}$ and $A^*(E\times_TU_N)=A^*(E\times_TV)/c_N(V)$, where $c_N$ is the $N-$th Chern class relative to $V$ (which is a $T$-characteristic class coming from $E$).

Fix now some $J\subseteq\{1,\dots, N\}$ with $|J|=p$ and let $\mathfrak{k}$ be the field of definition of $T_{\text{Orb}(J)}$ (which, we recall, is the subfield of $\overline{k}$ corresponding to the stabilizer of $J$ with respect to the action on $\{1,\dots, N\}$ of the Galois group of $k$).

We have $W_p=\underset{|J|=p}{\coprod} T_{\text{Orb}(J)}$, and by the localization exact sequence
\[
A^*(E\times_TW_p)\xrightarrow{i_*} A^*(E\times_TU_{p+1})\rightarrow A^*(E\times_TU_p)\rightarrow 0
\]
we see that we just need to prove that for each $J$ the image of $A^*(E\times_TT_{\text{Orb}(J)})$ in $A^*(E\times_TU_{p+1})$ lies in the image of $\mathcal{I}$.

Now observe that if $V_J$ is the standard representation over $\mathfrak k$ of $T_{\text{Orb}(J)}$ we have an open immersion $j:E\times_TT_{\text{Orb}(J)}\hookrightarrow E_\mathfrak{k}\times_{T_\mathfrak{k}}V_J$. 

Moreover we can see $T_{\text{Orb}(J)}$ as an irreducible component of $T_{\text{Orb}(J)}\times_k\text{Spec}(\mathfrak k)$, since $T_{\text{Orb}(J)}$ is defined over $\mathfrak k$; so we have a closed immersion $i'$ given by
\[
E\times_TT_{\text{Orb}(J)}\hookrightarrow E_{\mathfrak k}\times_{T_\mathfrak{k}}(T_{\text{Orb}(J)}\times_k\text{Spec}(\mathfrak k))\hookrightarrow (E\times_{T}U_{p+1})\times_k \text{Spec}(\mathfrak k).
\]
In particular we can factorize the closed immersion $i:E\times_TT_{\text{Orb}(J)}\hookrightarrow E\times_TU_{p+1}$ as the composition 
\begin{eqnarray}
E\times_TT_{\text{Orb}(J)}\overset{i'}{\hookrightarrow} (E\times_{T}U_{p+1})\times_k \text{Spec}(\mathfrak k)\overset{\pi^\mathfrak{k}}{\twoheadrightarrow} E\times_T U_{p+1}
\end{eqnarray}

Now, there is a fiber diagram
\[
\begin{tikzcd}
E\times_TT_{\text{Orb}(J)}\arrow[hookrightarrow]{r}{i'}\arrow[hookrightarrow]{d}{j}&(E\times_{T}U_{p+1})\times_k \text{Spec}(\mathfrak k)\arrow[hookrightarrow]{d}{j'}\\
E_\mathfrak{k}\times_{T_\mathfrak{k}}V_J\arrow[hookrightarrow]{r}{k}&E_{\mathfrak k}\times_{T_\mathfrak{k}}V_\mathfrak{k}
\end{tikzcd}
\]
where $k$ is the closed immersion induced by $V_J\hookrightarrow V_\mathfrak{k}$ and the vertical rows are open immersions.

Clearly we have that $j^*:A^*(M_\mathfrak{k})\simeq A^*(E_\mathfrak{k}\times_{T_\mathfrak{k}}V_J)\rightarrow A^*(E\times_TT_{\text{Orb}(J)})$ is surjective, and that $j'^*$ is also surjective (it is true for every open immersion). 

But $i'_*j^*=j'^*k_*$, so that if $x'\in A^*(M_\mathfrak{k})$ lifts $x\in A^*(E\times_TT_{\text{Orb}(J)})$, we can apply the self-intersection formula to see that $x'\cdot x_J$ lifts the push-forward of $x$ in $A^*((E\times_{T}U_{p+1})\times_k \text{Spec}(\mathfrak k))$. Here $x_J$ denotes the characteristic class in $A^*(M_\mathfrak{k})$ corresponding to $x_J\in A^*(BT_{\mathfrak k})$ (with the notation of the preceding section). Indeed, the conormal bundle of $E_\mathfrak{k}\times_{T_\mathfrak{k}}V_J$ in $E_\mathfrak{k}\times_{T_\mathfrak{k}}V_{\mathfrak k}$ is induced exactly by the standard representation $V_{J^c}$ of $T_{\text{Orb}(J^c)}$. This is easily seen recalling that $T_{\text{Orb}(J)}\times T_{\text{Orb}(J^c)}=T_{\mathfrak k}$ (\text{Lemma 1.2}). Its higher Chern class is exactly represented by $x_J$.\\

From the factorization $(1)$, it follows that $\pi^\mathfrak{k}_*(x'\cdot x_J)$ lifts $i_*(x)$ in $A^*(M)$. This clearly lies in $\mathcal{I}$, so we have our claim.
\end{proof}

\section{Some remarks on other nice classes of tori}

We can observe that the permutation tori $R_{E/k}\mathbb G_m$ are special groups, that is every $T-$torsor is Zariski-locally trivial: indeed by the main result (Theorem $1.1$) of \cite{Re} it suffices to check that for every field extension $k\subseteq K$ with Galois group $\Gamma:=\text{Gal}(\overline K:K)$ we have $0=H^1(K, K\otimes_k R_{E/k}\mathbb G_m)=H^1(K, R_{E\otimes_kK/K}\mathbb G_m)=H^1(\Gamma, (E\otimes_k\overline K)^*)$; but $E\otimes_k\overline K$ as a $\Gamma-$module splits as a product of $\overline K$, so the latter group is $0$ by Hilbert $90$ theorem.\\

Let us briefly see what is the structure of $A^*_T$ for such tori.

The ring $A^*_{T_{k'}}$ is easily seen (since $T_{k'}$ is a split torus) to be isomorphic to the polynomial ring $\mathbb{Z}[x_1,\dots, x_N]$, where each of the $x_i$'s corresponds to a character of $T_{k'}$. 

Indeed Totaro defines $A^*_G$ as the limit of $A^*(U/G)$ where $U$ ranges over open sets of $G-$representations on which $G$ acts freely; for the rest of the paper we will refer to these as \emph{approximating varieties for $BG$}.

The space $BT_{k'}$ is then the "limit" of the quotients $(V_{k'}-\{0\})^{\otimes n}/T_{k'}\simeq (\mathbb{P}^{N-1}_{k'})^{\otimes n}$, whose Chow rings are generated exactly by the (image of) elements $x_1,\dots, x_N$.

In particular we can identify the Chow ring $A^*_{T_{k'}}$ with the symmetric ring over the character module $\widehat T_{k'}$, which we will denote $S(\widehat T_{k'})$. We will call it \emph{the ring of characters} of the torus.\\

Now, $BT\simeq BT_{k'}/\Gamma$ is approximated by the spaces $(\mathbb{P}^{N-1}_{k'})^{\otimes n}/\Gamma$, where $\Gamma$ acts tautologically on $\text{Spec }k'$ and permutes the factors $\mathbb{P}^{N-1}_k$ according to its action on $V$.

In his article \cite{Kar}, Karpenko considers the approximating varieties for the classifying space of a permutation torus. In particular he takes the tori $T=R_{k'/k}\mathbb G_{m,k}$ obtained by Weil restriction from the multiplicative group. These approximating varieties are thus equal to the Weil restrictions $R_{k'/k}\mathbb P^n$ (where as before $k'$ is a splitting field for the permutation torus $T$). 

Then he shows how to calculate the class of these varieties in the ring of Chow motives  (for a definition and the main properties of Grothendieck's theory of motives see \cite{Dem}).

In the ring of motives, setting $\text{pt}=(\text{Spec }L,\text{Id})$ we have $\mathbb P^i_L=(\text{pt},0)\oplus(\text{pt},1)\dots\oplus (\text{pt},i)$ (\cite{Dem}, ex. 6); in particular we have a decomposition for $\mathbb P^i_L\times\dots \mathbb P^i_L$ in factors of the form $(\text{pt},k)$.\\

Karpenko determines the quotient of this by $\Gamma$; indeed, each factor of the above product is determined by a map $p:\{1,\dots,n\}\rightarrow\{0,1,\dots, i\}$, and $\Gamma$ permutes the factors corresponding to the $\Gamma-$orbits in the set of such maps. If $S$ is an orbit, for each $p\in S$ the quantity $|p|=\sum_{k=1}^n p(k)$ is the same, and we call this sum $|S|$. Karpenko's result states that

\begin{proposition}
For each orbit $S$ let $\mathfrak \Gamma\subseteq\Gamma$ be its stabilizer, and $L_S\subseteq L$ the corresponding subfield of $L$. The the motivic class of $R_{L/k}\mathbb P^i_L$ is equal to $\underset{S}{\bigoplus}(\emph{Spec }L_S,|S|)$.
\end{proposition}

By using this it is immediate to conclude that

\begin{corollary}
If $T$ is the Weil restriction $T=R_{E/k}\mathbb G_{m,k'}$ and $n=E:k$ then the map $A^*_T\rightarrow \mathbb Z[x_1,\dots, x_n]^\Gamma$ is an isomorphism.
\end{corollary}

Actually with the results of Section $1$ we could re-demonstrate this result.\\

Now, let us come back to special tori. Actually, we can characterize all the tori which are special (e.g. \cite{Hu}, Th. 5.1):

\begin{theorem}
An algebraic torus is special if and only if there exists a torus $T_1$ such that $T\times T_1$ is a permutation torus.
\end{theorem}

Note that the sufficiency is clear, since in this case for any field extension $K:k$ we have $H^1(K, T)\oplus H^1(K,T_1)=H^1(K, T\times T_1)=0$.

In this case, letting $P:=T\times T_1$, we have $BT\times BT_1=BP$ (we mean that, for each classifying bundles $U/T, U_1/T_1$, the product $U\times U_1/T\times T_1$ is an approximating variety of $BP$). In particular we have maps $i:BT\rightarrow BP$ and $p:BP\rightarrow BT$ such that $p\circ i=\text{Id}$, and $i^*\circ p^*=\text{Id}:A^*_T\rightarrow A^*_T$. We conclude that $p^*:A^*_T\rightarrow A^*_P$ is injective, so $A^*_T$ does not contain torsion classes since $A^*_P$ does not.\\

Now we want to show that the base change $A^*_T\rightarrow (A^*_{T_{k'}})^\Gamma$ (where $k'$ is a splitting field for $T$ and $\Gamma=\text{Gal}(k':k)$) is an isomorphism.

We immediately note that this is always true after tensoring with $\mathbb Q$ (this is actually true for any quotient by a finite group; see \cite{Fu}, 1.7.6).

Recall that if $T\times T_1=P$ is a permutation torus and $U/T,U_1/T_1$ are approximating varieties, then $(U\times U_1)/P$ is an approximating $P$-bundle. We would like to show that, passing to the limit, the image of $p^*:A^*_T\rightarrow A^*_P$ contains all the symmetric tensors coming from $S(\widehat T)\subset S(\widehat P)$.\\

We now use a result of Kimura (\cite{Ki}, Th. 2.3) on the action on the Chow rings of proper envelopes $\tilde X\xrightarrow{\pi} X$: in order to have $\alpha\in\pi^*(A^*(X))$ it is necessary and sufficient that $(\pi_1^*-\pi_2^*)(\alpha)=0$ in the following diagram
\[
\begin{tikzcd}
\tilde X\times_X \tilde X\arrow{r}{\pi_1}\arrow{d}{\pi_2}&\tilde X\arrow{d}\\
\tilde X\arrow{r}&X
\end{tikzcd}
\]

In the present case Kimura's criterion applies, when the base field is infinite, since $(U\times U_1)/P\simeq U/T\times U_1/T_1\rightarrow U/T$ is an envelope: indeed $U_1/T_1$ has a $k-$rational point (as $k$ is infinite), so that we can lift every subvariety of $U/T$; moreover we can choose $U/T$ such its base change to $\overline k$ is a product of projective spaces, hence proper.\\

Let $c$ be a symmetric tensor coming from $S(\widehat T)$: since the thesis holds rationally we know that $nc$ lies in the image of $p^*$ for some positive integer $n$, so from Kimura's criterion we have that $n(\pi_1^*-\pi_2^*)(c)=0$ in $BP\times_{BT}BP$ (here we mean that we take $(U\times U_1)/P\times_{U/T}(U\times U_1)/P$).

But $BP\times_{BT}BP=BP\times BT_1=B(T\times T_1)$ is the classifying space of a special torus, so its Chow ring is torsion free by the above observations. In particular for some $U,U_1$ we will have that $(\pi_1^*-\pi_2^*)(c)$ is not a torsion element in $(U\times U_1)/P\times_{U/T}(U\times U_1)/P$, so it must be zero. Then, again by Kimura's criterion we have that $c$ lies in the image of $p^*$. We conclude

\begin{proposition}
Let $T$ be a special torus over an infinite field with splitting field $k'$ and ${S(\widehat{T_{k'}})}$ be its ring of characters.

If $\Gamma=\text{\emph{Gal}}(k':k)$, then the map $A^*_T\rightarrow{S(\widehat{T_{k'}})}^\Gamma$ is an isomorphism. 
\end{proposition}

This also allows us to calculate the Chow ring of $BT$ when $T$ is an algebraic torus of the simplest kind apart from special tori, that is when there exists a resolution
\[
0\rightarrow T\rightarrow Q\rightarrow P\rightarrow 0
\]
with $P,Q$ special tori. Of course, without loss of generality we can assume that $P$ is a permutation torus: indeed if it is not the case then at least we can choose a special torus $P'$ such that $P\oplus P'$ is permutation, and then
\[
0\rightarrow T\rightarrow Q\oplus P'\rightarrow P\oplus P'\rightarrow 0
\]
is a resolution of the required kind.

Now, the point is that we have a map $BT\rightarrow BQ$ which is a $P-$torsor (to be precise: for every approximating torsor $U\rightarrow U/Q$ we have that $T$ acts freely on $U$ and $U/T\rightarrow U/Q$ is a $Q/T=P-$torsor).

In particular we can apply \textbf{Theorem 2.1} and \textbf{Corollary 3.2} to calculate $A^*_T$ in terms of $A^*_P,A^*_Q$. From what we know so far we can infer that $A^*_P=S(\widehat{P})^\Gamma$, where $\Gamma$ is the Galois group for some common splitting field over $k$ of $T,Q,P$, and that $A^*_Q\rightarrow A^*_T$ is surjective.

The ideal $\mathcal{I}=\ker(A^*_Q\rightarrow A^*_T)$ can be easily calculated. Indeed by the previous theorem we have that $\mathcal{I}$ contains, for every subgroup $\Gamma'<\Gamma$ corresponding to a field extension $L$, the image of $(A^*_{P_L})^+A^*_{Q_L}\rightarrow A^*_Q$ (where $M^+$ denotes the ideal of elements with positive degree in a graded ring). The push-forward from $A^*_{Q_L}$ to $A^*_Q$ is just $\text{Ind}^\Gamma_{\Gamma '}(-)$, where "$\text{Ind}^\Gamma_{\Gamma '}$" means the sum over representatives of the right cosets of $\Gamma/\Gamma '$; this is well-defined since $A^*_{Q_L}$ is inveriant under the action of $\Gamma '$. These elements are all the generators, and we can conclude the following

\begin{corollary}
If $T$ is an algebraic torus with a resolution $0\rightarrow T\rightarrow Q\rightarrow P\rightarrow 0$, $Q$ special and $P$ quasi-split, the ideal $\mathcal{I}=\ker(A^*_Q\twoheadrightarrow A^*_T)$ can be written as $\mathcal I=\underset{\Gamma '<\Gamma}{\sum}\text{\emph{Ind}}^\Gamma_{\Gamma '}(({S(\widehat{P})^+})^{\Gamma '}S(\widehat Q)^{\Gamma '})$.
\end{corollary}

This comes almost trivially from \textbf{Theorem 2.1}, since it implies that we have $\mathcal I=\underset{\Gamma '<\Gamma}{\sum}S(\widehat Q)^\Gamma\cdot\text{Ind}^\Gamma_{\Gamma '}(({S(\widehat{P})^+})^{\Gamma '}S(\widehat Q)^{\Gamma '})$. But for any $x\in S(\widehat Q)^\Gamma$ and any $y\in ({S(\widehat{P})^+})^{\Gamma '}S(\widehat Q)^{\Gamma '}$ we have $x\cdot \text{Ind}^\Gamma_{\Gamma '}(y)=\text{Ind}^\Gamma_{\Gamma '}(x\cdot y)$.\\

Let $\mathcal J:=S(\widehat{P})^+S(\widehat Q)$ be the kernel of the projection map $S(\widehat Q)\rightarrow S(\widehat T)$; then we have a commuting square
\[
\begin{tikzcd}
A^*_Q\arrow{r}{\simeq}\arrow[d, two heads]&S(\widehat Q)^\Gamma\arrow[d, two heads]\\
A^*_T\arrow{r}&S(\widehat T)^\Gamma
\end{tikzcd}
\] 
which gives us that $\text{coker}(A^*_T\rightarrow S(\widehat T))^\Gamma=\text{coker}(S(\widehat Q)^\Gamma\rightarrow S(\widehat T)^\Gamma)=H^1(\Gamma,\mathcal J)$.

The last step follows from the exact sequence
\[
S(\widehat Q)^\Gamma\rightarrow S(\widehat T)^\Gamma\rightarrow H^1(\Gamma,\mathcal J)\rightarrow H^1(\Gamma, S(\widehat Q))=0
\]
and the final equality is just Shapiro's lemma.

Finally $\ker(A^*_T\rightarrow S(\widehat T)^\Gamma)=\ker(A^*_Q/\mathcal I\rightarrow S(\widehat T)^\Gamma)=\mathcal J^\Gamma/\mathcal I$.\\

Let us summarize the above results:

\begin{proposition}
Let us be given a short exact sequence of tori $0\rightarrow T\rightarrow Q\rightarrow P\rightarrow 0$, whith $Q$ special and $P$ quasi-split. Then if $\phi$ is the map $A^*_T\rightarrow (A^*_{T_{k'}})^\Gamma$ and $\mathcal I,\mathcal J$ are defined as above we have that
\[
\ker(\phi)=\mathcal J^\Gamma/\mathcal I
\]
and
\[
\text{\emph{coker}}(\phi)=H^1(\Gamma,\mathcal I)
\]
\end{proposition}

By using these formulas we will be able to show that in some cases we have nonzero kernels and cokernels for the map $A^*_T\rightarrow S(\widehat T)^\Gamma$.

\section{Torsion and cotorsion classes}

By using our \text{Proposition 3.6} we will see how the kernel and cokernel of the map $A^*_T\rightarrow {A^*_{T_{\overline k}}}^\Gamma$ behave.\\

Special tori are not the largest class for which the base change map is an isomorphism. A very easy example is that of norm one tori, that is the closed subschemes $T=\ker N^{k'}_k\colon R_{k'/k}\mathbb G_{m,k'}\rightarrow\mathbb G_{m,k}$, whose $k-$rational points correspond to the elements in $k'^*$ with unitary norm over $k$. These are tori of the type considered at the end of the previous section, but let us nonetheless calculate directly their $A^*_T$.\\

These tori have character lattice $x_1\mathbb Z\oplus\dots\oplus x_n\mathbb Z/(x_1+\dots+x_n)$. Moreover, if $P=R_{k'/k}\mathbb G_{m,k'}$, then the norm map gives an isomorphism $P/T\simeq\mathbb G_{m,k}$; then $BT$ is a $P/T\simeq\mathbb G_{m,k}-$bundle over $BP$, so his Chow ring is equal to $A^*_P/c_1$, where $c_1$ is the first Chern class of this bundle. If $A^*_P=\mathbb Z[x_1,\dots,x_n]^\Gamma$ then this class is $x_1+\dots+ x_n$, so $A^*_T=\mathbb Z[x_1,\dots,x_n]^\Gamma/(x_1+\dots+ x_n)$.

The invariants of the character ring of $T$ are exactly $(\mathbb Z[x_1,\dots,x_n]/x_1+\dots+ x_n)^\Gamma\simeq \mathbb Z[x_1,\dots,x_n]^\Gamma/(x_1+\dots+ x_n)$, as we can see from the exact sequence
\[
\mathbb Z[x_1,\dots, x_n]^\Gamma\xrightarrow{\cdot (x_1+\dots+x_n)}\mathbb Z[x_1,\dots, x_n]^\Gamma\rightarrow \Bigl(\frac{\mathbb Z[x_1,\dots, x_n]}{x_1+\dots+x_n}\Bigr)^\Gamma\rightarrow H^1(\Gamma, \mathbb Z[x_1,\dots, x_n])=0
\]
where the last group is zero by Shapiro's lemma.\\

Now, by using our calculation of $A^*_T$ for a torus $T$ with a resolution by special tori $0\rightarrow T\rightarrow Q\rightarrow P\rightarrow 0$ we give some example of tori with nonzero kernel and cokernel for $A^*_T\rightarrow S(\widehat T)^\Gamma$.\\

Our first example is a torus $T$ such that the above map has a nonzero cokernel in degree $2$.

If $k,L$ are such that $\text{Gal}(L:k)=S_n$ (the symmetric group of rank $n>3$), we pick a torus with character lattice $\mathbb Z^n=\langle a_1,a_2,\dots,a_n\rangle$ equipped with an action of the Galois group $\Gamma=S_n$  so that $\sigma\in S_n$ sends each element of the basis $(a_1,\dots,a_n)$ to $\text{sgn}(\sigma)a_{\sigma(i)}$ (where $\sigma(-)$ denotes the tautological action on a $n$-element set). 

This has a resolution $0\rightarrow \mathbb Z^n\rightarrow \mathbb Z^{2n}\rightarrow\widehat T\rightarrow 0$; the middle term $\mathbb Z^{2n}$ is generated by elements $a_i^+, a_i^-$ for $i=1,\dots, n$ with a $S_n-$action $\sigma(a_i^\pm)=a_{\sigma(i)}^{\pm\text{sgn}(\sigma)}$, and the projection to $\widehat T$ is just $a_i^\pm\rightarrow \pm a_i$.

The kernel of this projection is the $n-$dimensional lattice generated by the $a_i^++a_i^-$, with the $S_n-$action inherited by its inclusion in $\mathbb Z^{2n}$. \\

This resolution corresponds to a reverse resolution $0\rightarrow T\rightarrow Q\rightarrow P\rightarrow 0$, where $P,Q$ are permutation tori of dimension $n, 2n$ respectively.\\

With the notations of the previous section we claim that $H^1(S_n, \mathcal J)$ is nonzero, and in particular the map $A^*_T\rightarrow S(\widehat T)^\Gamma$ has a nonzero cokernel.

This can be seen very easily at first glance: indeed the element $x=\underset{i<j}{\sum}a_ia_j\in S(\widehat T)^{S_n}$ cannot be lifted to any element of $S(\widehat Q)^{S_n}$; on the other hand $2x$ can be lifted to $\underset{i<j}{\sum}(a_i^+a_j^++a_i^-a_j^-)\in S(\widehat Q)^{S_n}$.

In order to see it, let us calculate $S(\widehat Q)_2^{S_n}$ (where the subscript labels the degree); the symmetric group permutes the monomials $a_i^\pm a_j^\pm$ which are a basis of $S(\widehat Q)_2$, so an invariant element containing a monomial $n_{ij}a_i^\pm a_j^\pm$ must also contain all the monomials $n_{ij}a_{\sigma(i)}^{\pm\text{sgn}(i)}a_{\sigma(j)}^{\pm\text{sgn}(\sigma)}$. However if $n>3$ then $S_n$ acts transitively on the ordered couples $(i,j)$, and moreover for every $(i,j)$ there is an odd permutation leaving the ordered pair fixed (just take a trasposition $(a,b)$ with $\{a,b\}$ disjoint from $\{i,j\}$). We conclude that all the invariants are of the form $\sum_i((a_i^\pm)^2+(a_i^\mp)^2)$ or $\sum_ia_i^+a_i^-$ or $\sum_{i<j}(a_i^{\epsilon_1}a_j^{\epsilon_2}+a_i^{-\epsilon_1}a_j^{-\epsilon_2})$ with $\epsilon_1,\epsilon_2\in\{-1,1\}$. Thus the image of $S(\widehat Q)_2^{S_n}$ in $S(\widehat T)_2^{S_n}$ consists of elements of the form $a\sum_ia_i^2+2b\sum_{i<j}a_ia_j$ ($a,b\in\mathbb Z$), and none of these can be $x$.\\

Now let us investigate more closely the problem of the existence of torsion classes in $A^*_T$. This is quite harder than the analogous problem on cotorsion, but we can nonetheless build an example.

Let us fix a field extension $L:k$ with Galois group $\Gamma=Q_8$, the group of quaternions. We recall that it is a group of order $8$ generated by elements $-1,i,j,k$ subject to the relations $(-1)^2=1, i^2=j^2=k^2=ijk=-1$. 

The group $Q_8$ has a natural action on $\mathbb Z^4$, and we can consider it as the character lattice of a $4-$dimensional torus $T$. There is a surjection $\mathbb Z^8\rightarrow \widehat T$, where $Q_8$ acts on $\mathbb Z^8$ as in its regular representation (or, more concretely, we identify a basis of $\mathbb Z^8$ with the eight elements of $Q_8$, who has a tautological action on them). We will call these elements $e,e',x,x',y,y',z,z'$ (corresponding respectively to $1,-1,i,-i,j,-j,k,-k$). This defines a permutation torus $Q$ and the kernel of $\widehat Q\rightarrow \widehat T$ is the sublattice of $\widehat Q$ generated by $e+e', x+x',y+y', z+z'$, on which $Q_8$ has a permutation action induced by $Q_8/\{1,-1\}\simeq {\mathbb Z/2\mathbb Z}\times {\mathbb Z/2\mathbb Z}$.  It gives a permutation torus $P$.\\

Let us consider the element $$w=\text{Orb}_{Q_8}(xyz)=xyz+x'y'z'+e'zy'+ez'y+z'e'x+zex'+yx'e'+y'xe$$ in $S(\widehat Q)^{Q_8}$, which lies in the kernel of $S(\widehat Q)\rightarrow S(\widehat T)$. We claim that this element does not lie in the ideal $\mathcal I$.

First of all we prove that $\mathcal I_3$ (the set of elements in degree $3$) is generated (additively) by $ \text{Ind}^{Q_8}_{\mathbb Z/2\mathbb Z}((S(\widehat P)^+)^{\mathbb Z/2\mathbb Z}S(\widehat Q)^{\mathbb Z/2\mathbb Z})_3$, where $\mathbb Z/2\mathbb Z=\{-1,1\}$. 

We begin by enumerating these elements. We note that $S(\widehat P)$ is already invariant by the action of $\{-1,1\}$. The elements of $S(\widehat Q)_1^{\mathbb Z/2\mathbb Z}$ are just those coming from $\widehat P$, while the elements of $S(\widehat Q)_2^{\mathbb Z/2\mathbb Z}$ are generated by $ee',xx',yy',zz'$ and the polynomials of the form $ab+a'b', ab'+a'b, a^2+a'^2$, where $a\neq b$ are in $\{e,i,j,k\}$. 

Thus there are elements of the form $\text{Ind}^{Q_8}_{\mathbb Z/2\mathbb Z}(S(\widehat P)_3)$ (first type) and elements of the form $\text{Ind}^{Q_8}_{\mathbb Z/2\mathbb Z}((a+a')y)$ where $a\in\{e,x,y,z\}$ and $y\in S(\widehat Q)_2^{\mathbb Z/2\mathbb Z}$ (second type). Actually, any element of the first type is also a sum of elements of the second type, hence we will restrict our attention to the latter.

Let $\mathcal L$ be this group of degree $3$ elements; we have to verify that $\mathcal I_3=\mathcal L$.

Let us check that these elements generate each $\text{Ind}^{Q_8}_{\Gamma '}((S(\widehat P)^+)^{\Gamma '}S(\widehat Q)^{\Gamma '})_3$ with $\Gamma '<Q_8$.\\

Let $\Gamma '=Q_8$. We note that all the elements of $(S(\widehat P)^+)^{Q_8}_{\leq 3}$ in degree $1,3$ are $\mathbb Z/2\mathbb Z$-induced elements from $S(\widehat P)^+$. Indeed the three elements of order two in $\mathbb Z/2\mathbb Z\times \mathbb Z/2\mathbb Z$ send $e+e'$ to $x+x',y+y', z+z'$ respectively, so that the full action on the basis $\widehat P$ is transitive; in particular an invariant monomial in $e+e',x+x',y+y',z+z'$ must contain all of them and thus must have degree at least four. Moreover, if a monomial is invariant under the action of one of the three generators, suppose without loss of generality that this generator swaps $e+e',x+x'$ (as well as $y+y',z+z'$): then the monomial must have factors of the form $(e+e')(x+x')$ and $(y+y')(z+z')$; the only possibility is that it has degree $2$.

Now, every element of the form $x\cdot y$ with  $x\in (S(\widehat P)^+)^{Q_8}_1, y\in S(\widehat Q)^{Q_8}_2$ or $x\in (S(\widehat P)^+)^{Q_8}_3, y\in\mathbb Z$ must be of the form $\text{Ind}^{Q_8}_{\mathbb Z/2\mathbb Z}(x'y)\in\mathcal L$ with $x'\in S(\widehat P)^+$, since $x$ is induced by some $x'$. 

Finally, an element of $S(\widehat Q)^{Q_8}_1$ must be a multiple of $e+e'+x+x'+y+y'+z+z'$, so every element of $(S(\widehat P)^+)^{Q_8}_2\cdot S(\widehat Q)^{Q_8}_1$ is a sum of elements of the previous form. This concludes the case $\Gamma '=Q_8$.\\

Suppose $\Gamma '=\mathbb Z/4\mathbb Z=\{1,-1,i,-i\}$ (or some other subgroup obtained swapping $i$ with $j$ or $k$, which can be checked exactly in the same way). 

This case is similar to the preceding one. Indeed in degree $1$ the $\Gamma '-$invariants of $S(\widehat P)^+$ are $e+e'+x+x'$ and $y+y'+z+z'$, while in degree $3$ the invariants are $(e+e')(x+x')(y+y'+z+z')$ and $(e+e'+x+x')(y+y')(z+z')$, which are all $\{-1,1\}$-induced from $S(\widehat P)^+$. So every element $x\cdot y$ with  $x\in (S(\widehat P)^+)^{\Gamma '}_1, y\in S(\widehat Q)^{\Gamma '}_2$ or $x\in (S(\widehat P)^+)^{\Gamma '}_3, y\in\mathbb Z$ is $\mathbb Z/2\mathbb Z$-induced from an element of $S(\widehat P)^+S(\widehat Q)^{\Gamma '}$. This implies that $\text{Ind}^{Q_8}_{\Gamma '}(xy)$ is $\mathbb Z/2\mathbb Z$-induced from $S(\widehat P)^+S(\widehat Q)^{\Gamma '}$.

Finally every element of $(S(\widehat P)^+)^{\Gamma '}_2\cdot S(\widehat Q)^{\Gamma '}_1$ is of the form above, since $S(\widehat Q)^{\Gamma '}_1$ is generated by $e+e'+x+x'$ and $y+y'+z+z'$.\\

If $\Gamma '=\{-1,1\}$ the claim is tautological, so we are left with checking the case $\Gamma '=\{1\}$. This is almost trivial: indeed $\text{Ind}^{\mathbb Z/2\mathbb Z}_{1}(S(\widehat P)^+S(\widehat Q))\subseteq S(\widehat P)^+S(\widehat Q)^{\mathbb Z/2\mathbb Z}$, since $S(\widehat P)^+$ is $\mathbb Z/2\mathbb Z$-invariant. In particular $$\text{Ind}^{Q_8}_1(S(\widehat P)^+S(\widehat Q))=\text{Ind}^{Q_8}_{\mathbb Z/2\mathbb Z}(\text{Ind}^{\mathbb Z/2\mathbb Z}_{1}(S(\widehat P)^+S(\widehat Q)))\subseteq \mathcal L.$$

Now we can finally prove that $w\notin \mathcal L$. Note that we can write $\mathcal L=A\oplus B$, where $A$ is generated by $\text{Ind}^{Q_8}_{\mathbb Z/2\mathbb Z}((a+a')y)$ with $a\in\{e,x,y,z\}$ and $y\in S(\widehat Q)_2^{\mathbb Z/2\mathbb Z}$ is an element such that either $a$ appears in some of its monomials or $y$ is of the form $bb', b^2+b'^2$; then we have $B$, generated by $\text{Ind}^{Q_8}_{\mathbb Z/2\mathbb Z}((a+a')y)$, where $a\in\{e,x,y,z\}$ and $y\in S(\widehat Q)_2^{\mathbb Z/2\mathbb Z}$ is of the form $pq+p'q'$ or $pq'+p'q$, with $p,q\in \{e,x,y,z\}\backslash\{a\}$ distinct from each other. In other words $A$ is generated by the elements such that every monomial has factors intersecting at most two of the couples $\{e,e'\},\{x,x'\},\{y,y'\},\{z,z'\}$; similarly the elements of $B$ have monomials with factors in distinct couples.

Clearly if $w\in \mathcal L$ then $w\in B$. However we immediately see that the sum of the coefficients of each generator of $B$ is $16$, and so the sum of the coefficients of $w$ should be a multiple of $16$. This is obviously not true, since the sum is $8$. This gives us our sought torsion class.

Let us finally note that $2w\in\mathcal L$, so that the class has order $2$. Indeed 

\begin{eqnarray}
2w&=&2\cdot \text{Ind}^{Q_8}_{\mathbb Z/2\mathbb Z}(xyz+x'y'z')\nonumber\\
&=&2\cdot\text{Ind}^{Q_8}_{\mathbb Z/2\mathbb Z}(yz(x+x')-x'z(y+y')+x'y'(z+z'))\nonumber\\
&=&\text{Ind}^{Q_8}_{\mathbb Z/2\mathbb Z}((x+x')(yz+y'z')-(y+y')(x'z+xz')+(z+z')(xy+x'y'))\nonumber
\end{eqnarray}

\end{document}